\documentclass[titlepage,11pt]{article}
\usepackage{latexsym}
\usepackage{amsfonts}
\usepackage{amsmath}
\usepackage{amsthm}
\newcommand\blackslug{\hbox{\hskip 1pt \vrule width 4pt height 8pt depth 1.5pt
        \hskip 1pt}}
\newcommand\bbox{\hfill \quad \blackslug \bigbreak}
\def\cupcup{\cup\cdots\cup}

\def\LL{,\ldots,}
\newcommand{\vare}{\varepsilon}

\title{Strengthening R\"odl's theorem}
\author{Maria Chudnovsky\thanks{This material is based upon work supported in part by the U. S. Army
Research Office under grant   number W911NF-16-1-0404, and supported by  NSF grant DMS 1763817.}\\
Princeton University, Princeton, NJ 08544
\\
\\
Alex Scott\thanks{Research supported by EPSRC grant EP/V007327/1.}\\
Mathematical Institute, University of Oxford, Oxford OX2 6GG, UK
\\
\\
Paul Seymour\thanks{Supported by AFOSR grant
A9550-19-1-0187, and by NSF grant  DMS-1800053.}\\
Princeton University, Princeton, NJ 08544
\\
\\
Sophie Spirkl\thanks{We acknowledge the support of the Natural Sciences and Engineering Research
Council of Canada (NSERC), [funding reference number RGPIN-2020-03912].
Cette recherche a \'et\'e financ\'ee par le Conseil de recherches en sciences
naturelles et en g\'enie du Canada (CRSNG), [num\'ero de r\'ef\'erence
RGPIN-2020-03912].}\\
University of Waterloo, Waterloo, Ontario N2L3G1, Canada}

\date{}
\swapnumbers
\newtheorem{thm}{Theorem}[section]
\newtheorem{lemma}[thm]{Lemma}

\newcommand{\Proof}{\noindent{\bf Proof.}\ \ }

\begin{document}
\maketitle
\begin{abstract}
What can be said about the structure of graphs that do not contain an induced copy of some graph $H$?  R\"odl showed in the 1980s that every $H$-free graph has large parts that are very sparse or very dense.  More precisely, let us say that
a graph $F$ on $n$ vertices is {\em $\vare$-restricted} if either $F$ or its complement has maximum degree at most $\vare n$.
R\"odl proved that for every graph $H$, and every $\vare>0$, every $H$-free graph $G$ has a linear-sized set of vertices inducing an $\vare$-restricted graph. 
We strengthen R\"odl's result as follows: for every graph $H$, and all $\vare>0$, 
every $H$-free graph can be partitioned into a bounded number of subsets inducing $\vare$-restricted graphs.
\end{abstract}

\section{Introduction}

What can be said about the structure of graphs that do not contain an induced copy of some graph $H$?
In the 1980s, R\"odl \cite{rodl} showed that every $H$-free graph has large parts that are very sparse or very dense.

To say that more precisely we need some definitions.
Graphs in this paper are finite and without loops or parallel edges. If $G$ is a graph and $X\subseteq V(G)$, we denote the subgraph
of $G$ induced on $X$ by $G[X]$, and $\overline{G}$ denotes the complement graph of $G$. If $G,H$ are graphs, we say that
$G$ is {\em $H$-free} if no induced subgraph of $G$ is isomorphic to $H$.
For a graph $G$,  let us say $X\subseteq V(G)$ is {\em weakly $\vare$-restricted}
if one of $G[X], \overline{G}[X]$ has at most $\vare|X|^2$ edges; and that $X$  is {\em $\vare$-restricted} if one of the graphs $G[X]$, $\overline{G}[X]$ 
has maximum degree at most $\vare|X|$.   R\"odl~\cite{rodl} proved the following:
\begin{thm}\label{rodlthm2}
For every graph $H$, and all $\vare>0$, there exists $\delta>0$ such that for every $H$-free graph $G$,
there is a weakly $\vare$-restricted set $X\subseteq V(G)$
with $|X|\ge \delta|G|$.
\end{thm}
Every $\vare$-restricted set is weakly $\vare/2$-restricted, and every weakly $\vare/2$-restricted set has a 
subset of at least half its size that is $\vare$-restricted.  Thus an equivalent version of R\"odl's theorem is the following:
\begin{thm}\label{rodlthm}
For every graph $H$, and all $\vare>0$, there exists $\delta>0$ such that for every $H$-free graph $G$,
there is an $\vare$-restricted set $X\subseteq V(G)$
with $|X|\ge \delta|G|$.
\end{thm}
R\"odl's theorem is an easy consequence of Szemer\'edi's regularity lemma, and has proved extremely useful. For example, it is now a standard tool in approaching the Erd\H{o}s-Hajnal conjecture (see for instance the breakthrough paper~\cite{lagoutte}, where it was crucial, and much subsequent work).
A proof of \ref{rodlthm} not using the regularity lemma (and consequently with much better constants) was given by Fox and Sudakov~\cite{fox}.

In this paper, we are concerned with {\em partitions} of $H$-free graphs such that {\em every} vertex class is either sparse or dense.  It is easy to prove that $H$-free graphs can be partitioned into a bounded number of weakly $\vare$-restricted subsets:
\begin{thm}\label{mainthmother}
For every graph $H$, and all $\vare>0$,
there is an integer $N$ such that for every $H$-free graph $G$, there is a partition of $V(G)$
into at most $N$ weakly $\vare$-restricted subsets.
\end{thm}
This can be shown by applying \ref{rodlthm2} repeatedly to partition most of the vertices into weakly 
$\vare/2$-restricted subsets, and then adding the remaining vertices into the largest set.  (See Nikiforov~\cite{nikiforov} for a strengthening of \ref{mainthmother}.)

But what about partitions into sets that satisfy the stronger property of being $\vare$-restricted?
This is much harder, and the main result of this paper is the following:
\begin{thm}\label{mainthm}
For every graph $H$, and all $\vare>0$,
there is an integer $N$ such that for every $H$-free graph $G$, there is a partition of $V(G)$
into at most $N$ $\vare$-restricted	 subsets.
\end{thm}
This is significantly stronger than \ref{mainthmother}.

Here is a third statement midway between the last two: that under the same hypotheses, $V(G)$ is the union of at most a bounded number of $\vare$-restricted subsets (not
necessarily pairwise disjoint). This variation does not seem to be easy, although it does not imply \ref{mainthm} as far as we know.

Some remarks: sets of cardinality at most two are always $\vare$-restricted, and for \ref{mainthm} 
it is sometimes necessary to use some $\vare$-restricted subsets of cardinality at most two, even in graphs $G$
with $|G|$ large.  For example, let $G$ be a star
$K_{1,n}$ with $n$ large, and let $\vare<1/3$: then every $\vare$-restricted subset containing the centre of the star has cardinality 
at most two.  (Note that this is not the case for \ref{mainthmother}; for example, a large star is already weakly $\vare$-restricted.)

Second, our proof of \ref{mainthm} (and the proof in~\cite{pure4} of \ref{countnonedges}, which
we will need to apply) does not use the regularity lemma. Thus we anticipate that the number $N$ in \ref{mainthm} is significantly smaller
(as a function of $1/\vare$) than numbers that are produced via the regularity lemma, but we have not made an estimate for it.

If $A,B\subseteq V(G)$ are disjoint, we say that
$B$ is {\em $\vare$-sparse to $A$} (in $G$) if every vertex in $B$ has at most $\vare|A|$ neighbours in $A$;
and $B$ is {\em $\vare$-dense
to $A$} if $B$ is $\vare$-sparse to $A$ in $\overline{G}$.
The method of proof of \ref{mainthm} is via the following statement:

\begin{thm}\label{mainlemma}
For every graph $H$, and all $\vare, \eta, \theta >0$, there exists an integer $N$ such that, for every $H$-free graph $G$, there is
a partition of $V(G)$ into nonempty sets 
$$A_1\LL A_{m}, B_1\LL B_{m}, C_1\LL C_n,$$ 
where $m\le |H|^2$ and $n\le N$, such that:
\begin{itemize}
\item $A_1 \LL A_m$ and $C_1\LL C_n$ are $\vare$-restricted sets;
\item for $1\le i\le m$, $|B_i|\le  \eta|A_i|$;
\item for $1\le i\le m$, $B_i$ is either $\theta$-sparse or $\theta$-dense to $A_i$.
\end{itemize}
\end{thm}
We will prove this in section~\ref{sec:mainlemma}. In section \ref{sec:pathpartition} we prove another result, and combine these
two to deduce \ref{mainthm}.

Let us give an idea of how \ref{mainlemma} will be used to prove \ref{mainthm}. 
We show in section \ref{sec:pathpartition} that if we have a sequence of $\lceil 4/\vare\rceil+1$ disjoint 
subsets of $V(G)$, say $W_0\LL W_K$, with three properties (that each of $W_0\LL W_{K-1}$ is much bigger than $W_K$, that
each of $W_0\LL W_{K-1}$ is $\vare/4$-restricted, and that $W_{i+1}\cupcup W_K$ is very sparse or very dense to $W_i$ for each $i$)
then the union of all these sets can be partitioned into a small number of $\vare$-restricted sets. This result is rather easy.

The role of \ref{mainlemma} is to deduce something similar for successively shorter sequences of sets, say $W_0\LL W_{k}$, with 
similar hypotheses. (If we can prove it for sequences of length one, then the main result follows, using the one-term sequence $V(G)$.)
We will need to adjust the parameters of the sequence as $k$ becomes smaller; that is, adjust the value of $\vare'$
such that $W_0\LL W_{k-1}$ are $\vare'$-restricted, and adjust
the density or sparsity condition.

To apply \ref{mainlemma}, we will use induction on $\lceil 4/\vare\rceil-k$; let $W_0\LL W_{k}$ be the sequence we want to handle now. We apply
\ref{mainlemma} to $G[W_{k}]$. That partitions $W_k$ into a bounded number of $\vare$-sparse sets and the pairs $(A_s,B_s)\; (1\le s\le m)$
as in \ref{mainlemma}. We could apply induction to the sequence $W_0\LL W_{k-1}, A_s,B_s$
for each $s$, and deduce that the union of $W_0\LL W_{k-1}, A_s,B_s$ can be partitioned into a bounded number of $\vare$-restricted 
sets for each $s$; but that does not help, since altogether these sets will cover the vertices of $W_0\cupcup W_{k-1}$ $m$ times instead of once. Instead, we first
partition each of $W_0\LL W_{k-1}$ into $m$ large subsets, one for each of the pairs $(A_1,B_1)\LL (A_m,B_m)$ that we need to 
handle; and then apply
induction to the sequence $W_0'\LL W_{k-1}',A_s,B_s$ for each $s$, where $W_j'$ is the ``large'' subset of $W_j$
that corresponds to $(A_s,B_s)$. The reason we partition each $W_j$ into ``large'' subsets, is that then the new sequences
will have parameters not much worse than the sequence $W_0\LL W_{k-1}, A_s,B_s$, and so we can apply the inductive hypothesis to them.

\section{Proving the main lemma}\label{sec:mainlemma}

In this section we prove \ref{mainlemma}. Let $A,B\subseteq V(G)$ be disjoint, and let $c,\vare>0$. 
We say that $(A,B)$ is {\em $(c,\vare)$-full} if
for all $A'\subseteq A$ with $|A'|\ge c|A|$ and $B'\subseteq B$ with $|B'|\ge c|B|$, the number of edges between $A',B'$ is at least
$\vare|A'|\cdot|B'|$. Similarly, $(A,B)$ is {\em $(c,\vare)$-empty} if it is $(c,\vare)$-full in the complement graph.
Thus if $(A,B)$ is $(c,\vare)$-full, and $A'\subseteq A$ and $B'\subseteq B$ with $|A'|/|A|, |B'|/|B|\ge c'> c$ then 
$(A', B')$ is $(c/c', \vare)$-full.

We need a version of a standard result called the ``embedding lemma'':

\begin{lemma}\label{transversal}
Let $G,H$ be graphs, let $0<\vare \le 1/2$, and let $A_v\;(v\in V(H))$ be pairwise disjoint nonempty subsets of $V(G)$, such that
for all distinct $u,v\in V(H)$, if $u,v$ are adjacent in $H$ then $(A_u, A_v)$ is $(\vare^{|H|},\vare)$-full, and if $u,v$ 
are nonadjacent then $(A_u, A_v)$ is $(\vare^{|H|},\vare)$-empty.
Then for each $v\in V(H)$ there exists $a_v\in A_v$ such that the map sending $v$ to $a_v$ for each $v\in V(H)$ is an isomorphism
from $H$ to an induced subgraph of $G$.
\end{lemma}
\Proof
We proceed by induction on $|H|$. If $|H|\le 1$ the result is true, so we assume $|H|>1$. Let $v\in V(H)$, and let $N,M$ be the sets
of neighbours of $v$ in $H$ and in $\overline{H}$ respectively. Let $c=\vare^{|H|}$. For each $u\in N$
there are fewer than $c|A_v|$ vertices in $A_v$ with fewer than $\vare|A_u|$ neighbours in $A_u$, since $(A_v, A_u)$ is
$(c,\vare)$-full; and similarly for each $u\in M$ there are fewer than $c|A_v|$ vertices in $A_v$ with fewer than $\vare|A_u|$ 
non-neighbours in $A_u$. Since $(|H|-1)c<1$ (because $\vare\le 1/2$), there exists $a_v\in A_v$ with at least  $\vare|A_u|$ neighbours in $A_u$ for each $u\in N$,
and at least $\vare|A_u|$ 
non-neighbours in $A_u$ for each $u\in M$. For each $u\in N$ let $B_u$ be the set of neighbours of $v$ in $A_u$, and for each $u\in M$
let $B_u$ be the set of non-neighbours of $v$ in $A_u$. Thus each $B_u\ne \emptyset$, since $|B_u|\ge \vare|A_u|$. 
Let $H'$ be obtained from $H$ by deleting $v$.

Thus for all distinct $u,w\in V(H')$,
if $u,w$ are adjacent then $(B_u, B_w)$ is $(c\vare^{-1}, \vare)$-full, and  if
$u,w$ are nonadjacent then $(B_u, B_w)$ is $(c\vare^{-1}, \vare)$-empty. From the inductive hypothesis, 
for each $u\in V(H')$ there exists $a_u\in B_u\subseteq A_u$ such that the map sending $u$ to $a_u$ for each $u\in V(H')$ is an isomorphism
from $H'$ to an induced subgraph of $G$. But then the theorem holds. This proves \ref{transversal}.~\bbox

The following is proved (without using the regularity lemma) in~\cite{pure4}, theorem 2.2:
\begin{lemma}\label{countnonedges}
For all $c,\vare,\tau>0$ with $\vare<\tau\le 8/9$, there exists $\gamma>0$ with the following property. Let $G$
be a bipartite graph with a bipartition $(A,B)$, with at least $\tau|A|\cdot|B|$  edges and with $A,B\ne \emptyset$. Then
there exist $A'\subseteq A$  and $B'\subseteq B$ with $|A'|/|A|, |B'|/|B|\ge \gamma$, such that
$(A',B')$ is $(c,\vare)$-full.
\end{lemma}

Now we are ready to prove \ref{mainlemma}, but first let us sketch its proof. We want to partition $V(G)$ into a bounded number of
$\vare$-restricted sets and ``$\theta$-restricted pairs'', by which we mean pairs of disjoint sets $(A,B)$ where $A$ is $\vare$-restricted, 
$B$ is much smaller
than $A$, and $B$ is either $\theta$-dense or $\theta$-sparse to $A$.
Let $V(H)=\{v_1\LL v_{|H|}\}$, and let us choose $t\le |H|$ maximum such that there is
a bounded-size collection $\mathcal{A}$ of $\vare$-restricted sets and $\theta$-restricted
pairs (the bound increasing with $t$), such that the set of vertices not in any member of $\mathcal{A}$
is partitioned into sets  $D_1\LL D_t, E$
satisfying three conditions:
\begin{itemize}
\item for all distinct $i,j$ with $1\le i,j\le t$, if $v_i,v_j$ are adjacent in $H$ then $(D_i,D_j)$ is $(x,y)$-full,
and  if $v_i,v_j$ are nonadjacent in $H$ then $(D_i,D_j)$ is $(x,y)$-empty, for some appropriate
 $x,y$;
\item $D_1\LL D_t$ are nonempty and $\vare'$-restricted where $\vare'$ is very small; and
\item $D_1\LL D_t$ are all many times bigger than the ``leftover'' set $E$.
\end{itemize}
We know that $t<|H|$, by \ref{transversal}, and now we will use $E$ to try to increase $t$ by 1. 
To build a new set $D_{t+1}$ within $E$, we like vertices
that are adjacent to at least a small fraction of the vertices in $D_i$ for the values of $i$ such that $v_i,v_{t+1}$ are 
adjacent in $H$, and are nonadjacent to at least a small fraction of the vertices in $D_i$ for the values of $i$ such that $v_i,v_{t+1}$ are 
nonadjacent (briefly, vertices that ``have the desired adjacency''). But the vertices that do not have the desired adjacency
are very sparse or very dense to some $D_i$,
and so we can remove them all from $E$ 
by adding a few more $\theta$-restricted pairs to $\mathcal{A}$. (We have to keep all the sets in $\mathcal{A}$ disjoint
from $D_1\LL D_t$, so we will have to shrink some of the sets a little, but that is straightfoward.)
So we can assume that every vertex in $E$ has the 
desired adjacency. Maybe $E=\emptyset$, but if so we have proved what we want, so we assume $E\ne \emptyset$, and this will lead to a contradiction.
By \ref{rodlthm} we can choose a linear subset $F_0$ of $E$ that is $\vare'$-restricted, where $\vare'$ is 
very small. By applying \ref{countnonedges} to each of the pairs $F_0,D_i$ in turn, we can choose a linear subset $D_{t+1}$
of $F_0$ that satisfies the first and second bullets above (changing $\vare'$ appropriately). 
We need to add a few sets to $\mathcal{A}$ to satisfy the third bullet.
There are two main issues to worry about. 
\begin{itemize}
\item 
First, when we applied \ref{countnonedges} to $F_0,D_i$, both $F_0, D_i$ might shrink by a constant factor, and we need to take care of the ``lost''
vertices, those that belong to the old $D_i$ but not the new one.
But the old $D_i$ was  $\vare'$-restricted where $\vare'$ is very small, and so we can arrange that the set of lost vertices is
$\vare$-restricted, by making sure it is not too small, and then we can add it to $\mathcal{A}$.
\item Second, we have to arrange 
that the new leftover set, $E'$ say, is small compared with $D_1\LL D_{t+1}$. The sets $D_1\LL D_t$ were shrunk in the process
of finding $D_{t+1}$; but they remain at least a constant factor of their original sizes, and so their sizes are at least some 
constant times $|E|$. And the same is true for $D_{t+1}$, since $D_{t+1}$ contains at least a linear fraction of $F_0$,
and $F_0$ contains a linear fraction of $E$.
But $E'$ is a subset of $E$, 
so, while  $E'$ might be bigger than some of $D_1\LL D_{t+1}$, its size is at most some large 
constant times the smallest of $D_1\LL D_{t+1}$; and therefore, by repeatedly applying \ref{rodlthm} and adding the sets we find to $\mathcal{A}$, we can reduce 
the size of $E'$ by any constant factor that we wish, and so bring its size down to what we need.
\end{itemize}

We restate \ref{mainlemma}:
\begin{thm}\label{mainlemma2}
For every graph $H$, and all $0<\vare, \eta,\theta<1$, there exists an integer $N$ such that, for every $H$-free graph $G$, there is
a partition of $V(G)$ into nonempty sets 
$$A_1\LL A_{m}, B_1\LL B_{m}, C_1\LL C_n,$$ 
where $m\le |H|^2$ and $n\le N$, such that:
\begin{itemize}
\item $A_1 \LL A_m$ and $C_1\LL C_n$ are $\vare$-restricted sets;
\item for $1\le i\le m$, $|B_i|\le  \eta|A_i|$;
\item for $1\le i\le m$, $B_i$ is either $\theta$-sparse or $\theta$-dense to $A_i$.
\end{itemize}
\end{thm}
\Proof
We may assume that $\vare,\eta,\theta<1/3$, by reducing them if necessary. For each $\vare'>0$, let $\delta_{\vare'}$ satisfy \ref{rodlthm} with $\vare,\delta$ replaced 
by $\vare', \delta_{\vare'}$.

Let $\vare_{|H|}=\min(\vare,(\theta/4)^{|H|})$. For $t=|H|-1, |H|-2\LL 0$ in turn:
\begin{itemize}
\item for $i = t, t-1\LL 0$ in turn, let $\Gamma_{t,i}=\gamma_{t,i+1}\Gamma_{t,i+1}$ (or $1$ if $i=t$);
and choose $\gamma_{t,i}$ 
such that \ref{countnonedges} holds, with $c,\vare,\tau, \gamma>0$ replaced by $\Gamma_{t,i}\vare_{t+1}/3, \theta/4,\theta/2, \gamma_{t,i}$
respectively (by decreasing $\gamma_{t,i}$ if necessary we may assume that $\gamma_{t,i}\le 1/3$ and $\gamma_{t,i}\le \gamma_{t,i+1}$);
\item let 
$\vare_t=\gamma_{t,0}c_{t+1}$.
\end{itemize}

For each $\gamma>0$, let $\phi(\gamma)$ be the smallest nonnegative integer that satisfies $(1-\delta_\vare)^{\phi(\gamma)}\le \gamma$.
For $0\le t\le |H|$, define 
$$\ell_t=\sum_{1\le i\le t}\phi(q_i),$$
where $p_i=\vare_{i}\Gamma_{i-1,0}$ and $q_i=\min(\gamma_{i-1,0},\frac12 \eta\delta_{p_i}\Gamma_{i-1,0})$.
Let $N= \ell_{|H|} + |H|(|H|+1)/2$; we claim that $N$ satisfies the theorem.

Let $G$ be $H$-free.
\\
\\
(1) {\em For all $\gamma$ with $0<\gamma<1$, and 
for every $X\subseteq V(G)$, there is a partition
of $X$ into at most $\phi(\gamma)+1$ sets, so that one of them has cardinality at most $\gamma|X|$ and the others are all $\vare$-restricted.}
\\
\\
Let $X\subseteq V(G)$. Choose an $\vare$-restricted set $A_1\subseteq X$ with $|A_1|\ge \delta_{\vare} |X|$; and inductively
for each $i>1$, choose an $\vare$-restricted set $A_i\subseteq X\setminus (A_1\cupcup A_{i-1})$ with 
$|A_i|\ge \delta_{\vare} |X\setminus (A_1\cupcup A_{i-1})|$.
It follows that $|X\setminus  (A_1\cupcup A_{i})|\le (1-\delta_\vare)^i|X|$ for each $i\ge 0$, and in particular when $i=\phi(\gamma)$.
This proves (1).

\bigskip

Let $V(H)$ have vertices $v_1\LL v_{|H|}$. For $0\le t\le |H|$,
we are interested in partitions of $V(G)$ into (possibly empty) sets $A_1\LL A_m, B_1\LL B_{m}, C_1\LL C_\ell$, $D_1\LL D_t$, and $E$, with the 
following properties:
\begin{itemize}
\item $m\le t(t-1)/2$ and $\ell\le \ell_t$;
\item $A_1\LL A_m, C_1\LL C_\ell$ and $D_1\LL D_t$ are all nonempty;
\item $A_1\LL A_m, C_1\LL C_{\ell}$ are $\vare$-restricted;
\item for $1\le i\le m$, $|B_i|\le \eta |A_i|$, and $B_i$ is either $\theta$-sparse or $\theta$-dense to $A_i$;
\item $D_1\LL D_t$ are $\vare_t$-restricted;
\item for $1\le i<j\le t$, if $v_i, v_j$ are adjacent in $H$ then $(D_i, D_j)$ is $(\vare_t, \theta/4)$-full, and if 
$v_i, v_j$ are nonadjacent then $(D_i, D_j)$ is $(\vare_t, \theta/4)$-empty;
\item $|E|\le (\eta/2) \min(|D_1|\LL |D_t|)$ if $t>0$.
\end{itemize}
Let us call such a thing a {\em partition of type $(m,\ell,t)$}.
To make clear which set plays which role in the partition, we will write them as:
$$(A_1,B_1)\LL (A_m, B_m)$$
$$C_1\LL C_\ell$$
$$D_1\LL D_t$$
$$E.$$
Choose such a partition,  of type $(m,\ell,t)$ say, with $t\le |H|$ maximum.
(This is possible, since $G$ admits a partition of type $(0,0,0)$, setting $E=V(G)$.)

Since $\vare_{|H|}\le  (\theta/4)^{|H|}$, and $D_1\LL D_t$ are nonempty, it follows from \ref{transversal} that $t\le |H|-1$.
Choose pairwise disjoint subsets $E_1\LL E_t$ of $E$ with maximal union, such that for $1\le i\le t$,
if $v_{t+1}, v_i$ are adjacent in $H$ then $E_i$ is $\theta/2$-sparse to $D_i$, and if $v_{t+1}, v_i$ are nonadjacent in $H$ then
$E_i$ is $\theta/2$-dense to $D_i$.
Let 
$E_0=E\setminus (E_1\cupcup E_t)$. Thus for $1\le i\le t$, $E_0$ is $(1-\theta/2)$-dense to 
$D_i$ if $v_{t+1}, v_i$ are adjacent in $H$,
and $E_0$ is $(1-\theta/2)$-sparse to $D_i$ if $v_{t+1}, v_i$ are nonadjacent.
Suppose, for a contradiction, that $E_0\ne \emptyset$.

We recall that $|E|\le \eta \min (|D_1|\LL |D_t|)$, and since $E\ne \emptyset$ (because $E_0\subseteq E$), it follows that $|D_i|\ge \eta^{-1}>1$ for $1\le i\le t$.
Thus $\lfloor |D_i|/2\rfloor\ge |D_i|/3$, for $1\le i\le t$.

Let $\vare'=\vare_{t+1}\Gamma_{t,0}$, and let $\delta'= \delta_{\vare'}$.
From \ref{rodlthm} there is an $\vare'$-restricted subset $F_0\subseteq E_0$ with $|F_0|\ge \delta'|E_0|$.
For $1\le i\le t$ define $F_i\subseteq F_{i-1}$ with 
$|F_i|\ge \gamma_{t,i}|F_{i-1}|$, and $H_i\subseteq E_i$ with $|D_i|/2\ge |H_i|\ge \gamma_{t,i}|D_i|$,
as follows. Let us assume that $v_{t+1}, v_i$ are adjacent (if they are non-adjacent, the construction is the same in the complement).
Thus $F_{i-1}$ is  $(1-\theta/2)$-dense to $D_i$. (We remark that this is a weak assertion: it means that each vertex in $F_{i-1}$
has at most $(1-\theta/2)|D_i|$ non-neighbours in $D_i$, but $\theta$ may be very small.) From the definition of $\gamma_{t,i}$, there exist $F_i\subseteq F_{i-1}$ and $H_i'\subseteq D_i$,
with $|F_{i}|\ge \gamma_{t,i}|F_{i-1}|$ and $|H_i'|\ge \gamma_{t,i}|D_i|$, such that $(F_i, H_i')$ is $(\Gamma_{t,i}\vare_{t+1}/3,\theta/4)$-full.
Let $H_i\subseteq H_i'$ of cardinality $\min(|H_i'|, \lfloor|D_i|/2\rfloor)$. Thus $|H_i|\ge \gamma_{t,i}|D_i|$, because either $|H_i|=|H_i'|\ge \gamma_{t,i}|D_i|$, or
$|H_i|=\lfloor|D_i|/2\rfloor\ge |D_i|/3\ge \gamma_{t,i}|D_i|$. Since $|H_i|\ge |H_i'|/3$, it follows that
$(F_i, H_i)$ is $(\Gamma_{t,i}\vare_{t+1},\theta/4)$-full.
This completes the inductive definition.

Thus, for $1\le i\le t$, $(F_i, H_i)$ is $(\Gamma_{t,i}\vare_{t+1},\theta/4)$-full if $v_i, v_{t+1}$ are adjacent, and 
$(\Gamma_{t,i}\vare_{t+1},\theta/4)$-empty if $v_i, v_{t+1}$ are non-adjacent.
Also, since $|F_i|\ge \gamma_{t,i}|F_{i-1}|$ for $1\le i\le t$, it follows that 
$|F_t|\ge \Gamma_{t,i}|F_i|.$
Consequently $(F_t,H_i)$ is $(\vare_{t+1},\theta/4)$-full if $v_i, v_{t+1}$ are adjacent, and  
$(\vare_{t+1},\theta/4)$-empty if $v_i, v_{t+1}$ are non-adjacent.

Now $|E|\le \eta \min (|D_1|\LL |D_t|)$. Let $\eta'=\min(\gamma_{t,0},\frac12 \eta\delta'\Gamma_{t,0})$. 
By (1) there exist pairwise disjoint, nonempty, $\vare$-restricted subsets $J_1\LL J_n$ of $E_0\setminus F_t$,
with $n\le \phi(\eta')$, such that their union
($J$ say) satisfies 
$$|E_0\setminus (F_t\cup J)|\le \eta'|E_0\setminus F_t|.$$
We claim that the sets 
$$(A_1,B_1)\LL (A_m,B_m),(D_1\setminus H_1, E_1)\LL (D_t\setminus H_t, E_t)$$
 $$C_1\LL C_{\ell}, J_1\LL J_n$$
$$H_1\LL H_t, F_t$$
$$E_0\setminus (F_t\cup J)$$
form a partition of $V(G)$ of type $(m+t,\ell+n,t+1)$.
To show this, we must check the following conditions: 

\begin{itemize}
\item Is it true that $m+t\le t(t+1)/2$ and $\ell+n\le \ell_{t+1}$? The first holds since $m\le t(t-1)/2$; and the second holds since 
$\ell+n\le \ell_t+\phi(\eta')=\ell_{t+1}$.
\item Is it true that $A_1\LL A_m, D_1\setminus H_1\LL D_t\setminus H_t,C_1\LL C_\ell, J_1\LL J_n, H_1\LL H_t$ and $F_t$ are all nonempty?
Certainly $A_1\LL A_m, C_1\LL C_{\ell}$ are nonempty from their definition, and so are $J_1\LL J_n$. For $1\le i\le t$, 
since $|E|\le \eta |D_i|$ and 
$E\ne \emptyset$, it follows that $|D_i|\ge 2$; and so $D_i\setminus H_i\ne \emptyset$, since $|H_i|\le |D_i|/2$.
Also $|H_i|=\min(|H_i'|, \lfloor|D_i|/2\rfloor)$, and $|H_i'|\ge \gamma_{t,i}|D_i|>0$, and $\lfloor|D_i|/2\rfloor>0$, so $H_i$ is
nonempty. Finally, $|F_t|\ge  \Gamma_{t,0}|F_0|$, and $|F_0|\ge \delta'|E_0|$, and $E_0\ne \emptyset$ by assumption; so $F_t\ne \emptyset$.
\item Is it true that $A_1\LL A_m, D_1\setminus H_1\LL D_t\setminus H_t,  C_1\LL C_{\ell}, J_1\LL J_n$ are $\vare$-restricted? $A_1\LL A_m$, 
$C_1\LL C_\ell$ and $J_1\LL J_n$ are $\vare$-restricted from their definition. For $1\le i\le t$, $D_i$ is $\vare_t$-restricted, and since
$|H_i|\le |D_i|/2$, it follows that $D_i\setminus H_i$ is $2\vare_t$-restricted and hence $\vare$-restricted.
\item Is it true that for $1\le i\le m$, $|B_i|\le \eta |A_i|$, and $B_i$ is either $\theta$-sparse or $\theta$-dense to $A_i$; and for $1\le i\le t$,
$|E_i|\le \eta |D_i\setminus H_i|$, and $E_i$ is either $\theta$-sparse or $\theta$-dense to $D_i\setminus H_i$? The first is true
from their definition. For the second, let $1\le i\le t$. Then 
$$|E_i|\le |E|\le (\eta/2)|D_i|\le \eta |D_i\setminus H_i|$$
since $|D_i\setminus H_i|\ge |D_i|/2$.  Also, $E_i$ is either $\theta/2$-sparse to $D_i$ (if $v_i, v_{t+1}$
are adjacent in $H$) or $\theta/2$-dense to $D_i$ (if $v_i, v_{t+1}$
are nonadjacent); and so $E_i$ is either $\theta$-sparse or $\theta$-dense to $D_i\setminus H_i$.
\item Is it true that  $H_1\LL H_t, F_t$ are $\vare_{t+1}$-restricted? For $1\le i\le t$, $D_i$ is $\vare_t$-restricted, and since $|H_i|\ge \gamma_{t,i}|D_i|\ge \gamma_{t,0}|D_i|$,
$H_i$ is $\vare_t/\gamma_{t,0}$-restricted and hence $\vare_{t+1}$-restricted. Also $F_0$ is $\vare'$-restricted, and $|F_t|\ge  \Gamma_{t,0}|F_0|$; and so 
$F_t$ is $\vare'/\Gamma_{t,0}$-restricted and hence $\vare_{t+1}$-restricted.
\item Is it true that for $1\le i<j\le t+1$, if $v_i, v_j$ are adjacent in $H$ then $(H_i, H_j)$ is $(\vare_{t+1}, \theta/4)$-full, and if
$v_i, v_j$ are nonadjacent then $(H_i, H_j)$ is $(\vare_{t+1}, \theta/4)$-empty? If $j=t+1$, we already saw that
$(F_t,H_i)$ is $(\vare_{t+1},\theta/4)$-full if $v_i, v_{t+1}$ are adjacent, and
$(\vare_{t+1},\theta/4)$-empty if $v_i, v_{t+1}$ are non-adjacent. So we may assume that $j\le t$.
Assume that $v_i, v_j$ are adjacent (the other case is similar). Then $(D_i, D_j)$ is $(\vare_{t}, \theta/4)$-full, and since
$|H_i|\ge \gamma_{t,0}|D_i|$ and $|H_j|\ge \gamma_{t,0}|D_j|$, and 
$\vare_t=\gamma_{t,0}\vare_{t+1}$, it follows that $(H_i, H_j)$ is $(\vare_{t+1}, \theta/4)$-full. 
\item Is it true that $|E_0\setminus (F_t\cup J)|\le (\eta/2) \min(|H_1|\LL |H_t|, |F_t|)$? For $1\le i\le t$, the choice of $J$
implies that 
$$|E_0\setminus (F_t\cup J)|\le \eta' |E_0\setminus F_t|\le \gamma_{t,0}|E_0|\le \gamma_{t,0}|E|\le (\gamma_{t,0}\eta/2)|D_i|\le (\eta/2)|H_i|.$$
Finally, to show that $|E_0\setminus (F_t\cup J)|\le (\eta/2) |F_t|$, observe that 
$$|E_0\setminus (F_t\cup J)|\le  \eta' |E_0|\le \eta'|E|\le \eta'|F_0|/\delta'\le (\eta/2)\Gamma_{t,0}|F_0|\le (\eta/2)|F_t|.$$
\end{itemize}

This proves that $G$ admits a partition of type $(m+t,\ell+n,t+1)$, contrary to the choice of $t$, and so 
completes the proof that $E_0=\emptyset$.

We claim that the pairs $(A_i,B_i)$ with $B_i\ne \emptyset$, the pairs $(D_i,E_i)$ with $E_i\ne \emptyset$, the sets $A_i$ 
(respectively, $D_i$) with $B_i$ (respectively, $E_i$) empty, and the sets $C_1\LL C_{\ell}$, satisfy the theorem. To show this, 
we may assume (by renumbering) that $B_1\LL B_r\ne \emptyset$, and $B_{r+1}\LL B_m=\emptyset$, and $E_1\LL E_s\ne \emptyset$, 
and $E_{s+1}\LL E_t=\emptyset$. We observe:
\begin{itemize}
\item The sets 
$$A_1\LL A_m, B_1\LL B_r, C_1\LL C_{\ell}, D_1\LL D_t, E_1\LL E_s$$
are pairwise disjoint and nonempty, and have union $V(G)$.
\item
The sets $A_1\LL A_m, C_1\LL C_{\ell}$ and $D_1\LL D_t$
are $\vare$-restricted (because each $D_i$ is $\vare_t$-restricted, and $\vare_t\le \vare$).
\item $|B_i|\le  \eta|A_i|$ for $1\le i\le r$, and 
$|E_i|\le |E|\le (\eta/2) \min(|D_1|\LL |D_t|)\le \eta|D_i|$ for $1\le i\le s$.
\item For $1\le i\le m$, $B_i$ is either $\theta$-sparse or $\theta$-dense to $A_i$, and for $1\le i\le s$,
$E_i$ is either $\theta/2$-sparse or $\theta/2$-dense to $D_i$, and hence either $\theta$-sparse or $\theta$-dense to $D_i$.
\item $r+s\le |H|^2$, since $r\le m\le t(t-1)/2$ and $s\le t$, and $t\le |H|$; and $\ell+(m-r)+(t-s)\le N$, since
$$\ell\le \ell_t\le \ell_{|H|}=N-|H|(|H|+1)/2$$
and 
$$m-r+t-s\le m+t\le t(t-1)/2+t \le |H|(|H|+1)/2.$$
\end{itemize}
This proves \ref{mainlemma2}.~\bbox

\section{Path-partitions}\label{sec:pathpartition}

We need the following two lemmas. For the first, see for example~\cite{bollobas}. 
\begin{lemma}\label{bollobas}
If $0\le k\le n$ are integers, then $\binom{n}{k} \le (en/k)^k$.
\end{lemma}
The second lemma is the following (logarithms in this paper are to base $e$):
\begin{lemma}\label{coveringlemma}
Let $\vare>0$ with $\vare\le 1/16$, and let $p\ge 0$ be an integer.
Let $G$ be a graph, and let $A,B$ be nonempty disjoint subsets of $V(G)$, such that 
$B$ is $\vare$-sparse to $A$, and 
$\log(2|B|)/\vare \le p\le |A|/12$.
Then there exists $P\subseteq A$ with $|P|=p$, such that $P$ is $2\vare$-sparse to $B$, and $B$ is $12\vare$-sparse to $P$.
\end{lemma}
\Proof
We may assume that some vertex in $B$ has a neighbour in $A$, because otherwise the result holds,
and since $\vare\le 1/16$ it follows that $|A|\ge 16$. Let $Q$ be the set of vertices in $A$ with fewer than
$2\vare |B|$ neighbours in $B$, and let $q=|Q|$.
There are at least $(|A|-q)(2\vare |B|)$ and at most $\vare |A|\cdot |B|$ edges between $A$ and $B$, and so
$q\ge |A|/2\ge 8$. Let $k=\lceil 12\vare p\rceil$.

Let $u_1\LL u_{2p}\in Q$, not necessarily all distinct. Let $y$ be the number of subsets of $Q$ of cardinality
$p$ that contain all of $u_1\LL u_{2p}$ (note that $p\le |A|/12\le q$);
and for each $v\in B$, let $z(v)$ be the number of subsets $I\subseteq \{1\LL 2p\}$
of cardinality $k$ such that $u_i$ is adjacent to $v$ for all $i\in I$ (note that $k= \lceil 12\vare p\rceil \le \lceil 2p\rceil = 2p$).
\\
\\
(1) {\em There is a choice
of $u_1\LL u_{2p}$ with $y=0$ and $z(v)=0$ for all $v\in B$.}
\\
\\
Choose $u_1\LL u_{2p}\in Q$ uniformly and
independently at random.
Let $\overline{y}$ be the expectation of $y$, and $\overline{z}(v)$ the expectation of each $z(v)$.
We will show that $\overline{y}<1/2$, and $\overline{z}(v)\le 1/(2|B|)$ for each $v\in B$, from which the claim follows.
First,
$$\overline{y}=\binom{q}{p}\left(\frac{p}{q}\right)^{2p}\le \left(\frac{ep}{q}\right)^p,$$
by \ref{bollobas}. 
Since $p\le |A|/12$ and $q\ge |A|/2$, it follows that $ep/q< 1/2$, and so
$\overline{y}< 1/2$.

For $v\in B$, since $v$ has at most $\vare |A|\le 2\vare|Q|$ neighbours in $Q$, it follows that
$$\overline{z}(v)\le \binom{2p}{k}(2\vare)^{k}\le 
\left(\frac{2ep}{k}\right)^k(2\vare)^{k}= \left(\frac{4e\vare p}{k}\right)^k\le \left(\frac{e}{3}\right)^k\le \left(\frac{e}{3}\right)^{12\vare p}$$
from \ref{bollobas}, and since $k\ge 12\vare p$ and $e<3$.
From the hypothesis, 
$\log(2|B|) \le \vare p\le 12\vare p\log(3/e)$, and so $\left(e/3\right)^{12\vare p}\le 1/(2|B|)$.
Hence $\overline{z}(v)\le 1/(2|B|)$,
and so the sum of $\overline{y}$ and all the $\overline{z}(v)\;(v\in B)$ is less than one.
This proves (1).

\bigskip

Choose $u_1\LL u_{2p}$ as in (1). Since $y=0$ it follows that $|\{u_1\LL u_{2p}\}|\ge p$; choose $P\subseteq \{u_1\LL u_{2p}\}$
with $|P|=p$.
Each vertex in $P$ has at most $2\vare |B|$ neighbours in $B$, since
$P\subseteq Q$; and each $v\in B$ has at most $12\vare p$ neighbours in $P$, since $z(v)=0$. This proves \ref{coveringlemma}.~\bbox

Let $G$ be a graph, let $k\ge 0$ be an integer, and let $\vare>0$. 
A {\em $(k,\vare)$-path-partition} of $G$ is a sequence $(W_0, W_1\LL W_{k})$ of subsets of $V(G)$, pairwise 
disjoint and with union $V(G)$, such that for $0\le i\le k-1$:
\begin{itemize}
\item $W_i$ is $\vare$-restricted;
\item $|W_k|\le |W_i|/12$;
\item $W_{i+1}\cupcup W_k$ is either $\vare/12$-sparse or $\vare/12$-dense to $W_i$.
\end{itemize}
If we are trying to partition $V(G)$ into $\vare$-restricted sets, and $G$ admits a  $(k,\vare)$-path-partition, 
then all but one of 
its sets
are $\vare$-restricted; the difficulty lies in handling the final set $W_k$.

\begin{thm}\label{coverpath}
Let $0<\vare\le 1/3$, and
let $G$ be a graph admitting a $(k,\vare/4)$-path-partition, where $k=\lceil 4/\vare\rceil$.
Then $V(G)$ can be partitioned into at most 
$2400\vare^{-2}$
$\vare$-restricted subsets.
\end{thm}
\Proof Let $(W_0\LL W_{k})$ be a  $(k,\vare/4)$-path-partition of $G$, let $p=|W_{k}|$, and $\vare'=\vare/48$. 
\\
\\
(1) {\em We may assume that $\log(2kp)\le \vare' p$.}
\\
\\
Suppose not; then $\log(2kp)> \vare p/48$, and since $k\le 4/\vare+1\le 13/(3\vare)$, it follows that
$$26p/(3\vare)\ge 2kp>e^{\vare p/48}\ge (\vare p/48)^3/6,$$
(because $e^x\ge x^3/3!$ for all $x>0$). We deduce that
$p^2\le 52\cdot 48^3/\vare^4$, and so $p\le 2398.5/\vare^2$.
Since $k\le 13/(3\vare)\le 1.5/\vare^2$,
the theorem holds, because
$V(G)$ is the union of $W_0\LL W_{k-1}$ and the $p$ singletons $\{v\}\;(v\in W_{k})$. This proves (1).
\\
\\
(2) {\em For $0\le i\le k$, there exists $C_i\subseteq W_i$ with $|C_i|=p$, such that for $0\le i\le k-1$, either
\begin{itemize}
\item $C_i$ is $2\vare'$-sparse to $C_{i+1}\cupcup C_{k}$, and $C_{i+1}\cupcup C_{k}$ is $12\vare'$-sparse 
to $C_i$, or 
\item $C_i$ is $2\vare'$-dense to $C_{i+1}\cupcup C_{k}$, and $C_{i+1}\cupcup C_{k}$ is $12\vare'$-dense
to $C_i$.
\end{itemize}
}
The choice of $C_i$ is inductive, as follows: let $C_{k}=W_{k}$, and now suppose that $0\le i\le k-1$, and  
$C_{i+1}\LL C_{k}$ are defined.
Let $B=C_{i+1}\cupcup C_{k}$. Thus $|B|=(k-i)p$ and $B$ is either $\vare'$-sparse or $\vare'$-dense to $W_i$ (because 
$(W_0\LL W_{k})$ is a $(k,12\vare')$-path-partition).
Moreover, 
$p=|W_k|\le |W_i|/12$.
Suppose first that $B$ is $\vare'$-sparse to $W_i$.
By (1), $\log(2|B|)\le \log(2kp)\le \vare' p$. 
By \ref{coveringlemma}, taking $A=W_i$, and replacing $\vare$ by $\vare'$, we deduce that there exists 
$C_i\subseteq W_i$ with $|C_i|=p$, such that $C_i$ is $2\vare'$-sparse to $B$, and $B$ is $12\vare'$-sparse to $C_i$.
Similarly, if $B$ is $\vare'$-dense to $W_i$, then \ref{coveringlemma} applied in $\overline{G}$ implies that
there exists 
$C_i\subseteq W_i$ with $|C_i|=p$, such that $C_i$ is $2\vare'$-dense to $B$, and $B$ is $12\vare'$-dense to $C_i$.
In either case, this completes the inductive definition of $C_0\LL C_{k}$, and so proves (2).

\bigskip

Now for $0\le i\le k-1$, either $C_i$ is $2\vare'$-sparse to $C_{i+1}\cupcup C_{k}$, or $2\vare'$-dense to $C_{i+1}\cupcup C_{k}$;
choose $I\subseteq \{0\LL k-1\}$ with $|I|\ge k/2$ such that either $C_i$ is $2\vare'$-sparse to $C_{i+1}\cupcup C_{k}$ for all $i\in I$, or
$C_i$ is $2\vare'$-dense to $C_{i+1}\cupcup C_{k}$ for all $i\in I$. Let $C=\bigcup_{i\in I\cup \{k\}}C_i$.
\\
\\
(3) {\em $C$ is $\vare$-restricted.}
\\
\\
To see this, suppose first that $C_i$ is $2\vare'$-sparse to $C_{i+1}\cupcup C_{k}$ for all $i\in I$.
Let $v\in C_j$ where $j\in I\cup \{k\}$, and let $I_1=\{i\in I: i<j\}$, and $I_2=\{i\in I\cup \{k\}: i>j\}$. 
Since $C_j$ is $2\vare'$-sparse to $C_{j+1}\cupcup C_{k}$, it follows that $v$ has at most $2\vare' p(k-j)\le \vare p(k-j)/4$ 
neighbours in $C_{j+1}\cupcup C_{k}$ (and hence at most the same number
in $\bigcup_{i\in I_2}C_i$). For each $i\in I_1$, since $C_{i+1}\cupcup C_{k}$ (and hence $C_j$) is $12\vare'$-sparse to $C_i$, it follows that
$v$ has at most $12\vare'p=\vare p/4$
neighbours in $C_i$; and therefore $v$
has at most $\vare pj/4$ neighbours in $\bigcup_{i\in I_1}C_i$. Since $v$ has at most $p$ neighbours in $C_j$, 
it follows
that $v$ has at most 
$$\vare p(k-j)/4 + \vare pj/4+p=  \vare p k/4 +p\le  \vare pk/2\le \vare |C|$$
neighbours in $C$ (here we use that $k\ge 4/\vare$ and $|C|\ge pk/2$),
and so $C$ is $\vare$-restricted. If $C_i$ is $2\vare'$-dense to $C_{i+1}\cupcup C_{k}$ for all $i\in I$, we use the same 
argument in the complement. 
This proves (3).

\bigskip
For each $i\in I$, since 
$|C_i|=|W_{k}|\le 3|W_i|/4$
and $W_i$ is $\vare/4$-restricted, it follows that $W_i\setminus C_i$ is $\vare$-restricted.
But then $V(G)$ admits a partition into the sets $W_i\;(i\in \{0\LL k-1\}\setminus I)$, the sets
$W_i\setminus C_i\;(i\in I)$, and $C$, and these sets are all $\vare$-restricted. This is a total of
$k+1\le 4/\vare+2$ sets, and $4/\vare+2\le 5/\vare\le 5/(3\vare^2)$. This proves \ref{coverpath}.~\bbox

Next we combine \ref{mainlemma} and \ref{coverpath} to prove an analogue of \ref{coverpath} for shorter and shorter sequences, and 
hence to prove \ref{mainthm}. We will show the following (the $-1$ at the end is for inductive purposes):
\begin{thm}\label{coversmallpath}
Let $H$ be a graph, and let $h=|H|^2$. Let $0<\vare\le 1/3$, and let $K=\lceil 4/\vare\rceil$. 
Let $N$ be as in \ref{mainlemma}, 
with $\vare,\eta,\theta$ replaced by $\vare/(4(2h)^K),1/(3h),\vare/(48(2h)^K)$ respectively.
Let $0\le k\le K$, and let 
$G$ be an $H$-free graph admitting a $(k, (2h)^{k-K}\vare/4)$-path-partition $(W_0\LL W_k)$.
Then $V(G)$ can be partitioned into at most $h^{K-k}(2400/\vare^2+N)-1$
$\vare$-restricted subsets.
\end{thm}
\Proof We may assume that $|H|\ge 2$ and so $h\ge 4$.
We proceed by induction on $K-k$. If $K-k=0$ then the result follows from \ref{coverpath}, so we assume that $k<K$, and the result
holds for $k+1$. 
By 
\ref{mainlemma}, there is
a partition of $W_k$ into 
nonempty sets 
$$A_1\LL A_{m}, B_1\LL B_{m}, C_1\LL C_{n},$$ 
where $m\le h$ and $n\le N$, such that: 
\begin{itemize}
\item $A_1 \LL A_{m}$ and $C_1\LL C_{n}$ are $\vare/(4(2h)^K)$-restricted sets;
\item for $1\le s\le m$, $|B_s|\le  |A_s|/(3h)$;
\item for $1\le s\le m$, $B_s$ is either $\vare/(48(2h)^K)$-sparse or $\vare/(48(2h)^K)$-dense to $A_s$. 
\end{itemize}
Let $X$ be the union of the sets $C_1\LL C_{n}$.
Since each of these sets is  $\vare/(4(2h)^K)$-restricted and hence $\vare$-restricted, it follows that $X$ can be partitioned into at most
$N$ $\vare$-restricted sets.
If $m=0$, the theorem holds, since $V(G)$ is the union of the sets $C_1\LL C_{n}$ and $W_0\LL W_{k-1}$, and $n\le N$ and $k\le 4/\vare$;
so we assume that $m>0$. Consequently $|W_k|\ge |A_1|\ge 3 h|B_1|\ge 3h$, and for $0\le i<k$, $|W_i|\ge 12 |W_k|\ge 36h$;
and hence $|W_i|/(2h)\ge 1$.
It follows that $\lceil |W_i|/(2h)\rceil\le  |W_i|/h$, and therefore there are $h$ pairwise disjoint subsets of $W_i$,
each of cardinality at least $|W_i|/(2h)$. Consequently we may choose subsets $W_i^1\LL W_i^m$ of $W_i$, pairwise disjoint and with union
$W_i$, and each of cardinality at least $|W_i|/(2h)$.
\\
\\
(1) {\em 
For $1\le s\le m$, $(W^s_0\LL W^s_{k-1}, A_s, B_s)$ is a $(k+1,(2h)^{k+1-K}\vare/4)$-path-partition
of $G[V_s]$, where $V_s=W^s_0\cupcup W^s_{k-1} \cup A_s\cup B_s$.}
\\
\\
To see this, we must show that 
\begin{itemize}
\item $A_s$ is $(2h)^{k+1-K}\vare/4$-restricted;
\item $|A_s|\ge 12 |B_s|$;
\item $B_s$ is either $(2h)^{k+1-K}\vare/48$-sparse or $(2h)^{k+1-K}\vare/48$-dense to $A_s$;
\end{itemize}
and also that 
for $0\le i\le k-1$:
\begin{itemize}
\item $W^s_i$ is $(2h)^{k+1-K}\vare/4$-restricted;
\item $|W^s_i|\ge 12 |B_s|$;
\item $W^s_{i+1}\cupcup W^s_{k-1}\cup A_s\cup B_s$ is either $(2h)^{k+1-K}\vare/48$-sparse or $(2h)^{k+1-K}\vare/48$-dense to $W^s_i$.
\end{itemize}

The first three statements are immediate from the definition of the pair $(A_s,B_s)$.
For the last three, 
let $0\le i\le k-1$. It follows that $W_{i}$ is $(2h)^{k-K}\vare/4$-restricted, and since
$|W^s_{i}|\ge |W_{i}|/(2h)$,
we deduce that  $W^s_{i}$ is $(2h)^{k+1-K}\vare/4$-restricted.

To show that $|W^s_{i}|\ge 12|B_s|$, observe that $|W_{i}|\ge 12 |W_k|\ge 12|A_s| \ge 36h|B_s|$,
and so 
$|W^s_i|\ge |W_i|/(2h)\ge 18|A_s|$.

Finally, to show that $W^s_{i+1}\cupcup W^s_{k-1}\cup A_s\cup B_s$ is either $(2h)^{k+1-K}\vare/48$-sparse or 
$(2h)^{k+1-K}\vare/48$-dense to $W^s_i$,
observe that, since  $(W_0\LL W_k)$ is a $(k, (2h)^{k-K}\vare/4)$-path-partition, it follows that
$W_{i+1}\cupcup W_k$ is either $(2h)^{k-K}\vare/48$-sparse or $(2h)^{k-K}\vare/48$-dense to $W_{i}$,
and hence so is 
$$W^s_{i+1}\cupcup W^s_{k-1}\cup A_s\cup B_s;$$ and therefore the latter is either  
$(2h)^{k+1-K}\vare/48$-sparse or $(2h)^{k+1-K}\vare/48$-dense to $W^s_{i}$, since $|W^s_i|\ge |W_i|/(2h)$. This proves (1).

\bigskip

From (1) and the inductive hypothesis, $V_s$ can be partitioned into at most 
$h^{K-k-1}(2400/\vare^2+N)-1$
$\vare$-restricted subsets, for $1\le s\le m$. Since the sets $V_1\LL V_m, C_1\LL C_n$ are pairwise disjoint and have union $V(G)$,
we deduce that $V(G)$ can be partitioned into at most 
$$h(h^{K-k-1}(2400/\vare^2+N)-1)+N\le h^{K-k}(2400/\vare^2+N)-1$$
$\vare$-restricted subsets. 
This proves \ref{coversmallpath}.~\bbox

To deduce \ref{mainthm}, we may assume that $\vare\le 1/3$, by reducing $\vare$ if necessary; then \ref{mainthm}
is immediate from \ref{coversmallpath} with $k=0$, applied to the $(0, (2h)^{-K}\vare/4)$-path-partition with one term $V(G)$.

Finally, we remark that Tung Nguyen~\cite{tung} has recently proved a strengthening of \ref{mainthm}, the following:
\begin{thm}\label{tungthm}
For every graph $H$, and all $\vare>0$,
there exist $C>0$ and an integer $N>0$  such that for every graph $G$, if $k$ denotes the number of distinct isomorphisms from
$H$ to induced subgraphs of $G$, then there exists $S\subseteq V(G)$ with 
$|S|\le C k^{1/|H|}$, such that there is a partition of $V(G\setminus S)$ into at most $N$ $\vare$-restricted sets.
\end{thm}
His proof is by a modification of the arguments of this paper.

\section*{Acknowledgement}
We are very grateful for a referee's report, which was immensely careful and helpful.


\begin{thebibliography}{99}
\bibitem{bollobas} B. Bollob\'as, {\em Modern Graph Theory}, page 216.
\bibitem{lagoutte} N. Bousquet, A. Lagoutte and S. Thomass\'e, ``The Erd\H{o}s-Hajnal conjecture for paths and
antipaths'', {\em J. Combinatorial Theory, Ser. B} {\bf 113} (2015), 261--264.
\bibitem{fox} J. Fox and B. Sudakov, ``Induced Ramsey-type theorems'', {\em Advances in Math.} {\bf 219} (2008), 1771--1800.
\bibitem{tung} T. Bguyen, ``A further extension of R\"odl's theorem'', manuscript August 2022.
\bibitem{nikiforov}  V. Nikiforov, ``Edge distribution of graphs with few copies of a given graph'', {\em Combin. Probab. Comput.} {\bf 15} (2006), 895--902.
\bibitem{pure4} A. Scott, P. Seymour and S. Spirkl, ``Pure pairs. IV. Trees in bipartite graphs'', submitted for publication, {\tt arXiv:2009.09426}.
\bibitem{rodl} V. R\"odl, ``On universality of graphs with uniformly distributed edges'',
{\em Discrete Math.} {\bf 59} (1986), 125--134.
\end{thebibliography}
\end{document}